\documentclass[preprints,article,accept,moreauthors,pdftex]{my-mdpi} 

\firstpage{1} 
\pubvolume{9}
\issuenum{19}
\articlenumber{2355}
\pubyear{2021}
\copyrightyear{2021}
\externaleditor{Academic Editor: Somayeh Mashayekhi} 
\datereceived{5 September 2021} 
\daterevised{16 September 2021} 
\dateaccepted{19 September 2021} 
\datepublished{22 September 2021} 
\hreflink{https://doi.org/10.3390/math9192355}
\doinum{10.3390/math9192355}

\pdfoutput=1

\Title{Optimal Control Problems Involving Combined Fractional Operators with General Analytic Kernels $^{\dagger}$}

\TitleCitation{Optimal Control Problems Involving Combined Fractional Operators with General Analytic Kernels}

\Author{Fa\"{\i}\c{c}al Nda\"{\i}rou $^{\ddagger}$\orcidA{} and Delfim F. M. Torres *$^{,\ddagger}$\orcidB{}}

\AuthorNames{Fa\"{\i}\c{c}al Nda\"{\i}rou and Delfim F. M. Torres}

\AuthorCitation{Nda\"{\i}rou, F.; Torres, D.F.M.}

\address[1]{Center for Research and Development in Mathematics and Applications (CIDMA), 
Department of Mathematics, University of Aveiro, 3810-193 Aveiro, Portugal; faical@ua.pt}

\corres{Correspondence: delfim@ua.pt; Tel.: +351-234-370-668}

\firstnote{This research is part of first author's Ph.D. project, 
which is carried out at the University of Aveiro under the
Doctoral Program in Applied Mathematics of Universities 
of Minho, Aveiro, and Porto (MAP-PDMA).} 

\secondnote{These authors contributed equally to this work.}	

\abstract{Fractional optimal control problems via a wide class 
of fractional operators with a general analytic kernel
are introduced. Necessary optimality conditions of 
Pontryagin type for the considered problem are obtained 
after proving a Gronwall type inequality as well as 
results on continuity and differentiability of 
perturbed trajectories. Moreover, a Mangasarian type
sufficient global optimality condition for the general 
analytic kernel fractional optimal control problem 
is proved. An illustrative example is discussed.}



\begin{document}


\section{Introduction}

Fractional Calculus, as a generalization of the traditional calculus 
through derivation and integration of an arbitrary order, is a rapidly 
growing field of mathematical research. Indeed, due to the existence 
of many different fractional operators in the literature, there is 
an interest in defining a more general class 
of fractional operators, which include existing operators as particular cases. 
This is important in the sense that with a general framework of operators 
it might be possible to establish a mathematical theory for this general formalism, 
rather than considering specific models with particular results. 
In~this direction, Fernandez, \"{O}zarslan and Baleanu proposed in 2019 
a fractional integral operator, based on a general analytic kernel, 
that includes a number of existing and known operators~\cite{arran}. 
Since this seminal work of 2019, several interesting results appeared, e.g.,
determination of source terms for fractional Rayleigh--Stokes equations 
with random data~\cite{MR4062054}, new analytic properties 
of tempered fractional calculus~\cite{MR3995283}, simulation 
of nonlinear dynamics with fractional neural networks arising 
in the modeling of cognitive decision making processes~\cite{MR4061392},
new numerical methods for variable order fractional nonlinear 
quadratic integro-differential equations~\cite{MR4083099},
and analysis of impulsive $\varphi$-Hilfer fractional differential 
equations~\cite{MR4141591}. Here, we investigate, for the first time 
in the literature, optimal control problems that involve a combined 
Caputo fractional derivative with a general analytic kernel 
in the sense of Fernandez, \"{O}zarslan and Baleanu.

The subject of combined fractional derivatives deals with the 
issue of combining the past and the future of the modelling process 
into one single operator. This is done by a convex combination 
of the left and right fractional derivatives. The idea was firstly 
introduced in~\cite{MR2846374} by Malinowska and Torres, 
following a previous idea of Klimek~\cite{MR1917624},
and then further investigated by the authors in~\cite{MR2870032,MR2944107,MR2861352}. 
See also~\cite{MR3199720,MR4113824,MR3787702} and references therein.
As mentioned in~\cite{MR2870032}, one advantage of combining fractional derivatives 
lie in the fact that they allow to describe a more general class of variational problems. 
Thus, it seems natural to consider optimal control problems that involve combined fractional 
derivatives with a general analytic kernel.

It should be mentioned that there is a rich literature on optimal control with 
fractional operators. Recent results include, for example:
(i) sensitivity properties of optimal control problems governed 
by nonlinear Hilfer fractional evolution inclusions 
in Hilbert spaces~\cite{MR4308222};
(ii) existence of a solution for a class of fractional delayed 
stochastic differential equations with non-instantaneous 
impulses and fractional Brownian motions~\cite{MR4305728};
(iii) an optimal control analysis of a fractional COVID-19 
epidemic model to minimize the infection 
and maximize susceptible individuals
under the Atangana--Baleanu fractional operator 
in the Caputo sense~\cite{MR4301833};
(iv) a Pontryagin maximum principle for optimal control problems 
with concentrated parameters for a degenerate differential equation 
with the Caputo operator~\cite{MR4287942}; etc. 
However, no available results on optimal control, 
with combined general analytic kernels, 
exist in the literature. For non-combined
general analytic kernels, optimal control results
are scarce and restricted to the recent 
publication~\cite{Ndairou:Torres:2021}. There, a weak version 
of Pontryagin's Maximum Principle for optimal control problems 
involving a general analytic kernel is given
but the emphasis is on the classical setting of the calculus
of variations (e.g., isoperimetric variational problems) and with
results valid only in the class of piece-wise continuous differentiable
state trajectories and piece-wise continuous controls~\cite{Ndairou:Torres:2021}.
In contrast, our current results are more general, being valid
in the class of absolutely continuous state trajectories
and $L^2$ controls. Moreover, results of~\cite{Ndairou:Torres:2021}
are valid only in the absence of constraints on the values
of the controls; that is, the controls take values in all the Euclidean space  
while here we are able to deal with more general and challenge
situations when the controls may take values
in any time-dependent close convex set of $L^2$.

The manuscript is organized as follows. Section~\ref{prelm} presents 
preliminary notions and results needed in the sequel, 
and follows the original results of the paper 
(Sections~\ref{basic:prpo} and \ref{sectionR}).
In Section~\ref{basic:prpo}, we prove two results
that are fundamental in the development of our work: 
a duality relation (Lemma~\ref{lemdual}) and
integration by parts formulas (Lemma~\ref{intpart}).
The main results appear then in Section~\ref{sectionR}, 
where we state and prove a Gronwall type inequality 
(Theorem~\ref{Gronwall}) and, as application to this inequality, 
we prove two results: a result on the continuity of solutions 
(see Lemma~\ref{cont} and Corollary~\ref{different}) 
and a necessary optimality condition of Pontryagin type to an optimal control problem 
with a general analytic kernel in the sense of Fernandez, \"{O}zarslan and Baleanu
(Theorem~\ref{theo}). We end Section~\ref{sectionR} by proving a sufficient condition 
for global optimality (Theorem~\ref{theosuff}). An example, illustrating 
the applicability of the obtained results, is given 
(see Examples~\ref{Pexemple} and \ref{suff:example}).
Finally, Section~\ref{sec:conc} give the main conclusions of the paper, 
including some possible future directions of research.


\section{Preliminaries}
\label{prelm}

In this section, we recall the definitions of fractional operators based 
on general analytic kernels and state some of their properties, 
relevant to our work.

\begin{Definition}[See~\cite{arran}]
\label{def1}
Let $[a, b]$ be a real interval, $\alpha$ be a real parameter in $[0, 1]$, 
$\beta $ be a complex parameter with non-negative real part, 
and $R$ be a positive number satisfying $R> (b-a)^{Re(\beta)}$. 
Let $A$ be a complex function, analytic on the disc $D(0, R)$, 
and defined on this disc by the locally uniformly convergent power series
\[
A(x)= \sum^{\infty}_{n=0} a_n x^n.
\] 
The left and right-sided fractional integrals with general analytic kernels  
of a locally integrable function $x: [a, b]\rightarrow \mathbb{R}$ 
(that is, $x\in L^1\left([a, b], \mathbb{R}\right)$) are defined by 
$$
^{A}I^{\alpha, \beta}_{a+} x(t)
:= \int^t_a(t-s)^{\alpha-1}A\left( (t-s)^{\beta}\right)x(s)ds
$$
and 
$$
^{A}I^{\alpha, \beta}_{b-} x(t)
:= \int^b_t(s-t)^{\alpha-1}A\left( (s-t)^{\beta}\right)x(s)ds,
$$
respectively.
\end{Definition}

\begin{Notation}
For any analytic function $A$ as in Definition~\ref{def1}, 
we define $A_{\Gamma}$ as
\[
A_{\Gamma}(x):= \sum^{\infty}_{n=0}a_n\Gamma(\beta n + \alpha)x^n.
\]
\end{Notation}

\begin{Lemma}[Series formula~\cite{arran}]
For any integrable function $x\in L^1\left([a, b], \mathbb{R}\right)$, 
the following uniformly convergent series formulas for 
$^{A}I^{\alpha, \beta}_{a+} x$ and for $^{A}I^{\alpha, \beta}_{b-} x$ 
as functions on $[a,b]$ hold:
$$
^{A}I^{\alpha, \beta}_{a+} x(t):= \sum^{\infty}_{n=0} a_n
\Gamma (\beta n + \alpha)^{RL}I^{\alpha + n\beta}_{a+}x(t)
$$
and 
$$ 
^{A}I^{\alpha, \beta}_{b-} x(t):=  \sum^{\infty}_{n=0} 
a_n\Gamma (\beta n + \alpha)^{RL}I^{\alpha + n\beta}_{b-}x(t),
$$
where $^{A}I^{\alpha + n\beta}_{a+} $ and $^{A}I^{\alpha + n\beta}_{b-}$ 
are, respectively, the left and right-sided Riemann--Liouville 
fractional integrals of order $\alpha + n\beta$.
\end{Lemma}

\begin{Lemma}[Theorem 2.5 of~\cite{arran}]
\label{el1operation}
With all notations as in Definition~\ref{def1}, we have a well bounded operator 
\[
^{A}I^{\alpha, \beta}_{a+}: L^1\left([a, b], 
\mathbb{R}\right) \rightarrow L^1\left([a, b], \mathbb{R}\right)
\]
for any fixed $\alpha$ and $\beta$ with $Re(\alpha), Re(\beta)\geqslant 0$. 
Moreover, the operator norm denoted by $\Vert \cdot \Vert$ is obtained as
\[
\left\Vert ^{A}I^{\alpha, \beta}_{a+} \right\Vert 
= \underset{f\in L^1\left([a, b], \mathbb{R} \right)}
\sup \frac{\left\Vert ^{A}I^{\alpha, \beta}_{a+} 
f \right\Vert_1}{\Vert f \Vert_1} = (b-a)^{\alpha}M, 
\quad M = \underset{|x|< (b-a)^{\alpha}}\sup A(x).
\]
In addition, similar results hold for the right-sided operator 
given in Definition~\ref{def1}, that is, the operator 
$^{A}I^{\alpha, \beta}_{b-}$ is bounded on 
$L^1\left([a, b], \mathbb{R} \right)$ 
with operator norm at most $(b-a)^{\alpha}M$.
\end{Lemma}

\begin{Lemma}[Semi group property~\cite{arran}]
\label{semigroup}
Let $a, b, A$ be as in Definition~\ref{def1}, 
and fix $\alpha_1, \alpha_2, \beta \in \mathbb{C}$ 
with non-negative real parts. The semigroup property 
\[
^{A}I^{\alpha_1, \beta}_{a+} \circ ^{A}I^{\alpha_2, \beta}_{a+} x(t) 
= ^{A}I^{\alpha_1 + \alpha_2, \beta}_{a+}x(t) 
=  ^{A}I^{\alpha_2 + \alpha_1, \beta}_{a+}x(t) 
= ^{A}I^{\alpha_2, \beta}_{a+} \circ ^{A}I^{\alpha_1, \beta}_{a+} x(t)
\]
is uniformly valid (regardless of $\alpha_1, \alpha_2, \beta $ and $x$) 
if and only if the condition 
\[
\sum_{m+n=k}a_n(\alpha_1, \beta)a_m(\alpha_2, \beta)
\Gamma(\alpha_1 + n\beta)\Gamma(\alpha_2 + n\beta)
= a_k(\alpha_1 + \alpha_2, \beta)\Gamma(\alpha_1 + \alpha_2+ k\beta)
\]
is satisfied for all non-negative integers $k$.
\end{Lemma}

A similar result to that of Lemma~\ref{semigroup} 
holds for the right-sided Riemann--Liouville fractional 
integral operator given in Definition~\ref{def1}.

Next, we give some recalls on fractional derivatives 
with general analytic kernel in the sense of Riemann--Liouville and Caputo.

Let $a, b, \alpha, \beta$ and $A$ be as in Definition~\ref{def1} 
and denote
$$
L^{\alpha, \beta}\left([a, b], \mathbb{R} \right)
:= \left\{ x\in L^1\left([a, b], \mathbb{R} \right)
:\, ^{\bar{A}}I^{1-\alpha, \beta}_{a+} x, \,
^{\bar{A}}I^{1-\alpha, \beta}_{b-}x 
\in AC\left([a, b], \mathbb{R} \right)\right\},
$$ 
where $AC\left([a, b], \mathbb{R} \right)$ represents 
the set of absolutely continuous functions on $[a, b]$.

\begin{Definition}[See~\cite{arran}]
The left and right Rieman--Liouville fractional derivatives 
with general analytic kernels, of a 
function $x\in L^{\alpha, \beta}\left([a, b], \mathbb{R} \right)$, are defined by
\[
^{A}_{RL}D^{\alpha, \beta}_{a+} x(t)
= \frac{d}{dt}\left( ^{\bar{A}}I^{1-\alpha, \beta}_{a+}x(t)\right) 
\quad \text{ and } \quad 
^{A}_{RL}D^{\alpha, \beta}_{b-} x(t)
= -\frac{d}{dt}\left(^{\bar{A}}I^{1-\alpha, \beta}_{b-}x(t)\right),
\]
where function $\bar{A}$ used on the right-hand side is an analytic function 
defined by 
$$
\bar{A}(x)=\sum_{n=0}^{\infty}\bar{a}_n x^n
$$ 
and such that $A_{\Gamma}\cdot \bar{A}_{\Gamma}=1$.
\end{Definition}

Now, let us denote by $AC^{\alpha, \beta}\left([a, b], \mathbb{R} \right)$  
the set of absolutely continuous functions that can be represented as 
\begin{equation}
\label{absolutefun}
x(t)= x(a) + ^{\bar{A}}I^{\alpha, \beta}_{a+}f(t) 
\quad \text{ and } \quad
x(t)= x(b) + ^{\bar{A}}I^{\alpha, \beta}_{b-}f(t),
\end{equation}
for some function $f\in   L^{\alpha, \beta}\left([a, b], \mathbb{R}\right)$.

\begin{Definition}
\label{defCaputo}
The left and right Caputo fractional derivatives with general analytic kernels, 
of a function $x \in AC^{\alpha, \beta}\left([a, b], \mathbb{R} \right)$, are defined by
$$
^{A}_{C}D^{\alpha, \beta}_{a+} x(t)
= \frac{d}{dt}\left[ ^{\bar{A}}I^{1-\alpha, \beta}_{a+}\left( x(t)-x(a) \right)\right] 
$$
and
$$ 
^{A}_{C}D^{\alpha, \beta}_{b-} x(t)
= -\frac{d}{dt}\left[^{\bar{A}}I^{1-\alpha, \beta}_{b-}\left(x(t)-x(b)\right)\right],
$$
where function $\bar{A}$ used on the right-hand side 
is an analytic function given by
$\bar{A}(x)=\sum_{n=0}^{\infty}\bar{a}_n x^n$ and such that 
$A_{\Gamma}\cdot \bar{A}_{\Gamma}=1$.
\end{Definition}

We would like to emphasize that, from \eqref{absolutefun} 
and by using the semi group property (Lemma~\ref{semigroup}), 
one has
\[
^{\bar{A}}I^{1-\alpha, \beta}_{a+}\left( x(t)-x(a) \right)
= ^{\bar{A}}I^{1, \beta}_{a+}f(t)\quad \text{ and } 
\quad ^{\bar{A}}I^{1-\alpha, \beta}_{b-}\left(x(t)-x(b)\right)
= ^{\bar{A}}I^{1, \beta}_{b-}f(t).
\]
Therefore, by Definition~\ref{defCaputo} and the fact that 
$f\in L^{\alpha, \beta}\left([a, b], \mathbb{R}\right)$, 
it is obvious that $^{A}_{C}D^{\alpha, \beta}_{a+}$ and 
$^{A}_{C}D^{\alpha, \beta}_{b-}$ belong to 
$L^1\left([a, b], \mathbb{R}\right)$.

\begin{Lemma}[See, e.g.,~\cite{MR2013877,MR0274683}]
\label{lemma:concave}
Let $h: \mathbb{R}^{n} \rightarrow \mathbb{R}$ be a continuously differentiable function. 
Then $h$ is a concave function if and only if it satisfies the so called gradient inequality:
\[
h(\theta_1)-h(\theta_2)\geq \nabla h(\theta_1)(\theta_1- \theta_2)
\]
for all $\theta_1, \theta_2 \in \mathbb{R}^n$.
\end{Lemma}


\section{Fundamental Properties}
\label{basic:prpo}

We prove rules of fractional integration by parts 
for the general analytic kernel operators. Firstly, we show 
a duality formula for the integral operator (Lemma~\ref{lemdual}). 
Then, we use the duality formula to prove the fractional integration 
by parts formulas for the general analytic kernel fractional operators
(Lemma~\ref{intpart}).

\begin{Lemma}[Duality relation]
\label{lemdual}
Let $\alpha$ and $\beta$ be as in Definition~\ref{def1}
with $|\alpha + \beta| \geq 1$.
For any functions $x(t)$ and $y(t)$, $t \in [a, b]$, 
with $x, y \in L^1\left([a, b], \mathbb{R}\right)$, 
the following duality relation holds:
\[
\int^b_a x(t)^{A}I^{\alpha, \beta}_{a+} y(t) dt 
= \int^b_a y(t)^{A}I^{\alpha, \beta}_{b-} x(t)dt.
\]
\end{Lemma}

\begin{proof}
By the series formula, we know that
\begin{equation}
\label{series}
\int^b_a x(t)^{A}I^{\alpha, \beta}_{a+} y(t) dt 
= \int^b_a x(t)\sum^{\infty}_{n=0}a_n 
\Gamma(\beta n + \alpha)^{RL}I^{\alpha + n\beta}_{a+} y(t)dt.
\end{equation}
Since the series in the right hand side of \eqref{series} 
is uniformly convergent, it follows that
\[
\int^b_a x(t)^{A}I^{\alpha, \beta}_{a+} y(t) dt 
= \sum^{\infty}_{n=0}a_n \Gamma(\beta n + \alpha) 
\int^b_a x(t) ^{RL}I^{\alpha + n\beta}_{a+} y(t)dt
\]
and, by the well-known duality of 
the Riemann--Liouville integral operators, 
which is valid under assumption $|\alpha + \beta| \geq 1$
(cf.~Lemma~2.7(a) of~\cite{MR2218073}), we have
\[
\int^b_a x(t) ^{RL}I^{\alpha + n\beta}_{a+} y(t)dt 
=  \int^b_a y(t) ^{RL}I^{\alpha + n\beta}_{b-} x(t)dt
\]
for any $n \in \mathbb{N}$, which leads to 
\[
\sum^{\infty}_{n=0}a_n \Gamma(\beta n + \alpha) 
\int^b_a x(t) ^{RL}I^{\alpha + n\beta}_{a+} y(t)dt 
= \sum^{\infty}_{n=0}a_n \Gamma(\beta n + \alpha)  
\int^b_a y(t) ^{RL}I^{\alpha + n\beta}_{b-} x(t)dt.
\]
Therefore, we obtain that
\[
\int^b_a x(t)^{A}I^{\alpha, \beta}_{a+} y(t) dt 
= \int^b_a y(t)^{A}I^{\alpha, \beta}_{b-} x(t)dt.
\]
This concludes the proof.
\end{proof}

\begin{Lemma}
\label{lemmaequiCaputo}
Let $x\in AC^{\alpha, \beta}\left([a, b], \mathbb{R} \right)$. 
The left and right sided Caputo fractional derivatives, 
as defined in Definition~\ref{defCaputo}, 
coincide with the following representation:
\[
^{A}_{C}D^{\alpha, \beta}_{a+} x(t)
=  ^{\bar{A}}I^{1-\alpha, \beta}_{a+}x^{'}(t)
\quad \text{ and } \quad 
^{A}_{C}D^{\alpha, \beta}_{b-} x(t)
= -\left( ^{\bar{A}}I^{1-\alpha, \beta}_{b-}x^{'}(t)\right).
\]
\end{Lemma}

\begin{proof}
We have, by Definition~\ref{defCaputo}, that
\[
^{A}_{C}D^{\alpha, \beta}_{a+} x(t)
= \frac{d}{dt}\left[ ^{\bar{A}}I^{1-\alpha, 
\beta}_{a+}\left( x(t)-x(a) \right)\right].
\]
Using the series formula (Lemma~\ref{series}), it follows that
\begin{align*}
^{A}_{C}D^{\alpha, \beta}_{a+} x(t) 
&= \frac{d}{dt}\left\{ \sum^{\infty}_{n=0}\bar{a}_n 
\Gamma(\beta n + 1-\alpha)^{RL}I^{\beta n + 1-\alpha }_{a+} 
\left( x(t)-x(a)  \right)\right\}\\
&= \frac{d}{dt}\left\{ \sum^{\infty}_{n=0}\bar{a}_n 
\Gamma(\beta n + 1-\alpha) \left[ \frac{1}{\Gamma(\beta n + 1-\alpha)} 
\int^t_a (t-s)^{\beta n -\alpha} \left( x(s)-x(a)  \right)ds \right] \right\}.
\end{align*}
Moreover, by the classical integration by parts formula, we have
\begin{align*}
\int^t_a (t-s)^{\beta n -\alpha} \left( x(s)-x(a)  \right)ds 
&= \left[-\frac{1}{\beta n + 1-\alpha}(t-s)^{\beta n+1 
-\alpha} \left( x(s)-x(a)  \right) \right]^t_a \\ 
&\quad + \int^t_a\frac{1}{\beta n + 1 -\alpha}(t-s)^{\beta n 
+ 1-\alpha}\frac{d}{ds}\left( x(s)-x(a)  \right)ds\\
&= \int^t_a\frac{1}{\beta n + 1 
-\alpha}(t-s)^{\beta n + 1-\alpha}x^{'}(s)ds.
\end{align*}
Therefore, 
\begin{multline}
^{A}_{C}D^{\alpha, \beta}_{a+} x(t) 
= \frac{d}{dt}\left\{ \sum^{\infty}_{n=0}\bar{a}_n 
\Gamma(\beta n + 1-\alpha) \right.\\
\left.\times \left[ \frac{1}{\Gamma(\beta n + 1-\alpha)} 
\int^t_a\frac{1}{\beta n + 1 -\alpha}(t-s)^{\beta n + 1-\alpha}x^{'}(s)ds \right] \right\}.
\end{multline}
Since the series is uniformly convergent, 
we can differentiate with respect to $t$ to obtain
\begin{align*}
^{A}_{C}D^{\alpha, \beta}_{a+} x(t) 
&= \sum^{\infty}_{n=0}\bar{a}_n 
\Gamma(\beta n + 1-\alpha) \left[ \frac{1}{\Gamma(\beta n + 1-\alpha)} 
\int^t_a(t-s)^{\beta n -\alpha}x^{'}(s)ds \right] \\
&= \sum^{\infty}_{n=0}\bar{a}_n 
\Gamma(\beta n + 1-\alpha) \left[ ^{RL}I^{\beta n + 1-\alpha}_{a+} x^{'}(t) \right]
= ^{\bar{A}}I^{1-\alpha, \beta}_{a+}x^{'}(t).
\end{align*}
In a similar way, it is possible to derive the right sided representation. 
\end{proof}

The next result is, as we shall see, an important tool for proving 
necessary optimality conditions to optimal control problems.
 
\begin{Lemma}[Integration by parts formula]
\label{intpart}
Let $\alpha$ and $\beta$ be as in Definition~\ref{def1}
with $|\alpha + \beta| \geq 1$,
$x \in  L^{\alpha, \beta}\left([a, b], \mathbb{R}\right)$ 
and $y \in  AC^{\alpha, \beta}\left([a, b], \mathbb{R} \right)$.  
Then, the following two formulas hold: 
\begin{equation}
\int^b_a x(t) ^{A}_{C}D^{\alpha, \beta}_{a+} y(t)dt 
= \left[y(t)^{\bar{A}}I^{1-\alpha, \beta}_{b-}x(t)\right]^b_a 
+ \int^b_a y(t) ^{A}_{RL}D^{\alpha, \beta}_{b-} x(t)dt
\end{equation}
and
\begin{equation}
\label{formulaInt}
\int^b_a x(t) ^{A}_{C}D^{\alpha, \beta}_{b^{-}} y(t)dt 
= \left[-y(t)^{\bar{A}}I^{1-\alpha, \beta}_{a^{+}}x(t)\right]^b_a 
+ \int^b_a y(t) ^{A}_{RL}D^{\alpha, \beta}_{a^{+}} x(t)dt.
\end{equation}
\end{Lemma}

\begin{proof}
By definition, 
\[
\int^b_a x(t) ^{A}_{C}D^{\alpha, \beta}_{a+} y(t)dt 
= \int^b_a x(t) ^{\bar{A}}I^{1-\alpha, \beta}_{a+}y'(t)dt
\]
and, by the duality formula of Lemma~\ref{lemdual}, it follows that 
\[
\int^b_a x(t) ^{\bar{A}}I^{1-\alpha, \beta}_{a+}y'(t)dt 
= \int^b_a y'(t) ^{\bar{A}}I^{1-\alpha, \beta}_{b-}x(t)dt.
\]
Using (standard) integration by parts, we obtain that
\[
\int^b_a y'(t) ^{\bar{A}}I^{1-\alpha, \beta}_{b-}x(t)dt 
= \left[y(t)^{\bar{A}}I^{1-\alpha, \beta}_{b-}x(t)\right]^b_a 
- \int^b_a y(t)\frac{d}{dt}\left( ^{\bar{A}}I^{1-\alpha, \beta}_{b-}x(t)\right) dt,
\]
which leads to the desired formula. 
The proof of \eqref{formulaInt} is similar.
\end{proof}


\section{Main Results} 
\label{sectionR}

In this section, following the definition given in~\cite{MR2846374}, 
we propose a combined fractional operator for a general analytic kernel  
and study the related optimal control problem.

\begin{Definition}
\label{def:combined}
Let $\gamma \in [0, 1]$. The combined fractional operator 
with general analytic kernel is defined by
\[
^{A}_{C}D^{\alpha, \beta, \gamma}_{a, b}
= \gamma ^{A}_{C}D^{\alpha, \beta}_{a^{+}} 
+ (1-\gamma) ^{A}_{C}D^{\alpha, \beta}_{b^{-}}. 
\]
\end{Definition}

Note that $^{A}_{C}D^{\alpha, \beta, 0}_{a, b}
= {^{A}_{C}D}^{\alpha, \beta}_{b^{-}}$ and 
$^{A}_{C}D^{\alpha, \beta, 1}_{a, b}
= {^{A}_{C}D}^{\alpha, \beta}_{a^{+}}$. The operator is obviously linear. 
Using Lemma~\ref{intpart}, with $x \in  L^{\alpha, \beta}\left([a, b], \mathbb{R}\right)$ 
and $y \in AC^{\alpha, \beta}\left([a, b], \mathbb{R} \right)$, 
we can easily establish the following integration by parts formula:
\begin{multline}
\label{combineIntpart}
\int^b_a x(t) ^{A}_{C}D^{\alpha, \beta, \gamma}_{a, b} y(t)dt 
= \gamma \left[y(t)^{\bar{A}}I^{1-\alpha, \beta}_{b^{-}}x(t)\right]^b_a\\
+ (1-\gamma)\left[-y(t)^{\bar{A}}I^{1-\alpha, \beta}_{a^{+}}x(t)\right]^b_a
+\int^b_a y(t) ^{A}_{RL}D^{\alpha, \beta, 1-\gamma}_{a,b} x(t)dt,  
\end{multline}
where $^{A}_{RL}D^{\alpha, \beta, \gamma}_{a,b} 
= \gamma ^{A}_{RL}D^{\alpha, \beta}_{a^{+}} 
+ (1-\gamma) ^{A}_{RL}D^{\alpha, \beta}_{b^{-}}$.

In the following subsection, we prove a new integral inequality 
of Gronwall type that will be useful to investigate continuity 
of solutions to our optimal control problem.


\subsection{Gronwall's Inequality}

Gronwall's inequality is an important integral inequality 
that is often used to prove qualitative and quantitative properties 
of solutions to differential equations. Very recently, there were several 
works devoted to this subject in the field of fractional calculus: see, 
e.g.,~\cite{MR4020061,MR3991009} and references therein. 
The next result is a new Gronwall type inequality for a fractional 
integral operator with a general analytic kernel. 

\begin{Theorem}[Gronwall's inequality]
\label{Gronwall}
Let $f(\cdot), u(\cdot) \in L^1\left([a, b], \mathbb{R}\right)$ be non-negative 
and $g(\cdot)$ be a non-negative monotonic increasing continuous 
function on $[a, b]$ satisfying 
$$
\underset{t \in [a, b]}\max g(t) < \frac{1}{(b-a)^{\alpha-1}M}
\quad
\text{ with } 
\quad
M=  \underset{|x|< (b-a)^{\alpha}}\sup A(x).
$$
If 
\begin{equation}
\label{inequality}
u(t)\leq f(t)+ g(t)\left(^AI^{\alpha, \beta, \gamma}_{a, b}u\right)(t), 
\quad t\in [a, b] \ a.e.,
\end{equation}
then for almost all $t \in [a, b]$ we have 
\[
u(t) \leq f(t)+ \sum_{k=1}^{\infty}[g(t)]^{k}\left[
\left(  ^AI^{\alpha, \beta, \gamma}_{a, b}\right)^{k}(f)\right](t),
\] 
where $^AI^{\alpha, \beta, \gamma}_{a, b}
= \gamma ^AI^{\alpha, \beta}_{a^{+}} + (1-\gamma)^AI^{\alpha, \beta}_{b^{-}}$ 
with $\gamma \in [0, 1]$ and we use the composition of operators' notation: 
\[
\left(  ^AI^{\alpha, \beta, \gamma}_{a, b}\right)^{k} 
= \circ^{k}_{i=1}\left( ^AI^{\alpha, \beta, \gamma}_{a, b}\right)_i.
\]
\end{Theorem}

\begin{proof}
Since the operator $^AI^{\alpha, \beta, \gamma}_{a, b}$ is a non-decreasing operator, 
as linear combination of non-decreasing operators, it follows that
\begin{equation*}
\begin{split}
\left( ^AI^{\alpha, \beta, \gamma}_{a, b}u\right)(t)
&\leq {^AI^{\alpha, \beta, \gamma}_{a, b}}\left(f(\cdot) 
+ g(t)\left(^AI^{\alpha, \beta, \gamma}_{a, b}u\right)\right)(t)\\
&= \left(^AI^{\alpha, \beta, \gamma}_{a, b}f\right)(t) 
+ g(t)\left[\left(^AI^{\alpha, \beta, \gamma}_{a, b}\right)^2(u)\right](t).
\end{split}
\end{equation*}
Now, we substitute this previous inequality into \eqref{inequality}, 
to obtain that
\[
u(t) \leq f(t) + g(t)\left(^AI^{\alpha, \beta, \gamma}_{a, b}f\right)(t) 
+ [g(t)]^2\left[\left(^AI^{\alpha, \beta, \gamma}_{a, b}\right)^2 (u)\right](t).
\]
Repeating this procedure up to $N$ times, we get
\[
u(t) \leq f(t) + \sum_{k=1}^{N-1}\left( g(t)\right)^k\left[
\left(^AI^{\alpha, \beta, \gamma}_{a, b} \right)^k(f)\right](t) 
+ \left(g(t)\right)^N\left[\left(^AI^{\alpha, \beta, \gamma}_{a, b}\right)^N (u)\right](t).
\]
Therefore, when $N\rightarrow \infty $, one has
\[
u(t) \leq f(t) + \sum_{k=1}^{\infty}\left( g(t)\right)^k\left[
\left(^AI^{\alpha, \beta, \gamma}_{a, b} \right)^k(f)\right](t) 
+ \lim_{N\rightarrow \infty}\left(g(t)\right)^N\left[
\left(^AI^{\alpha, \beta, \gamma}_{a, b}\right)^N (u)\right](t).
\]
To obtain the desired result, it remains to show that the series 
$$
\sum_{k=1}^{\infty}\left( 
g(t)\right)^k\left[\left(^AI^{\alpha, \beta, \gamma}_{a, b} \right)^k(f) \right](t) 
$$ 
converges and the limit 
$$
\underset{{N\rightarrow \infty}}\lim \left(g(t)\right)^N\left[
\left(^AI^{\alpha, \beta, \gamma}_{a, b}\right)^N (u)\right](t)
$$ 
is equal to zero. Let us study the composition of operators 
$\left(^AI^{\alpha, \beta, \gamma}_{a, b}\right)^k$. For this purpose, 
note that, by Lemma~\ref{el1operation}, the operators 
$^AI^{\alpha, \beta}_{a^{+}}$ and $^AI^{\alpha, \beta}_{b^{-}}$ 
are both bounded on $L^1\left( [a, b], \mathbb{R}\right)$. Therefore, 
because $L^1\left( [a, b], \mathbb{R}\right)$ is a norm vector space, 
we have that the linear combination  
$\gamma ^AI^{\alpha, \beta}_{a^{+}} + (1-\gamma)^AI^{\alpha, \beta}_{b^{-}}
= {^AI}^{\alpha, \beta, \gamma}_{a, b}$ 
is also a bounded operator on $L^1\left( [a, b], \mathbb{R}\right)$ 
having the same operator norm, precisely 
\[
\left\Vert ^AI^{\alpha, \beta, \gamma}_{a, b} f \right\Vert_1 
\leqslant (b-a)^{\alpha}M \Vert f \Vert_1, 
\quad f\in L^1\left( [a, b], \mathbb{R}\right).
\]
As a consequence, if $f \in L^1\left( [a, b], \mathbb{R}\right)$, 
then for any fixed integer $k$ we have that the composition 
$\left(^AI^{\alpha, \beta, \gamma}_{a, b}\right)^k (f) 
\in L^1\left( [a, b], \mathbb{R}\right) $ 
and is bounded in the sense that
\[
\left\| \left(^AI^{\alpha, \beta, \gamma}_{a, b}\right)^k 
f \right\|_1 \leqslant \left[(b-a)^{\alpha}M\right]^k \Vert f \Vert_1, 
\quad f\in L^1\left( [a, b], \mathbb{R}\right).
\]
Moreover, using the mean value theorem, we have that there exists 
$t\in [a, b]$ such that 
$$
\left| \left(^AI^{\alpha, \beta, \gamma}_{a, b}\right)^k (f)(t) \right| 
= \frac{1}{(b-a)}\left\| \left(^AI^{\alpha, \beta, \gamma}_{a, b}\right)^k 
f \right\|_1. 
$$
Hence, we obtain that 
\[
\left| \left(^AI^{\alpha, \beta, \gamma}_{a, b}\right)^k (f)(t) \right|  
\leqslant \frac{1}{(b-a)}\left[(b-a)^{\alpha}M\right]^k  \Vert f \Vert_1, 
\quad f\in L^1\left( [a, b], \mathbb{R}\right),
\]
and it follows that
\[
\left|\sum_{k=1}^{\infty}\left( g(t)\right)^k\left(
^AI^{\alpha, \beta, \gamma}_{a, b}\right)^k(f)(t)\right| 
\leq  \Vert f \Vert_1 \sum_{k=1}^{\infty}\left[(b-a)^{\alpha-1}MN\right]^k, 
\quad f\in L^1\left( [a, b], \mathbb{R}\right),
\]
where $N= \displaystyle{\underset{t \in [a, b]}\max g(t)}$.
Finally, the series converges since, by assumptions,  
$$
N< \frac{1}{(b-a)^{\alpha-1}M}. 
$$
Moreover, according to the necessary condition of convergence 
of an infinite series, one deduces that
$$
\underset{{k\rightarrow \infty}}\lim \left(g(t)\right)^k\, \left(
^AI^{\alpha, \beta, \gamma}_{a, b}\right)^k(f)(t)=0.
$$ 
This concludes the proof.
\end{proof}


\subsection{Applications}

In this subsection, we prove several important results: 
continuity of solutions to optimal control problems (Lemma~\ref{cont}), 
which is an application of our Gronwall's inequality (Theorem~\ref{Gronwall});
differentiability of the perturbed trajectories (Corollary~\ref{different});
and a necessary optimality condition of Pontryagin type to problem \eqref{bp}
(Theorem~\ref{theo}), which happens to be an application of the results 
on continuity, differentiability, and integration by parts. 
First of all, let us define the optimal control problem that we will be studying.

We consider an analytic kernel fractional optimal control problem, which consists 
in finding a control $u \in \Omega(t)\subseteq L^2\left([a, b], \mathbb{R}\right)$ 
and its corresponding state trajectory $x \in AC^{\alpha, \beta}$ solution to the 
following problem:
\begin{equation}
\label{bp}
\begin{gathered}
J[x(\cdot), u(\cdot)]= \int^{b}_{a} L\left(t, x(t), u(t)\right)dt \longrightarrow \max,
\\
^A_{C}D^{\alpha, \beta, \gamma}_{a, b}x(t)= f\left(t, x(t), u(t)\right), \quad t\in [a, b],
\quad |\alpha + \beta| \geq 1,\\
x(\cdot) \in AC^{\alpha, \beta}, 
\quad u(\cdot) \in \Omega(t) \subseteq L^2, \quad a.e. \quad t\in [a, b],\\
x(a)= x_a,
\end{gathered}
\end{equation}
where $\Omega(t)$ is a closed convex subset of $L^2$ and functions $L$ and $f$ 
are assumed to be continuously differentiable in all their three arguments, 
that is, $f \in C^1$ and $L\in C^1$. In~particular, $f$ is locally Lipschitz  
with Lipschitz constant $K$. By solution of the analytic kernel
fractional optimal control problem \eqref{bp}, we mean a pair 
$\left(x(\cdot), u(\cdot)\right) \in AC^{\alpha, \beta} \times \Omega(t)$
satisfying the control system $^A_{C}D^{\alpha, \beta, \gamma}_{a, b}x(t)
= f\left(t, x(t), u(t)\right)$, $t\in [a, b]$, the initial condition 
$x(a)= x_a$, and giving maximum value to functional $J$.
This solution is given by Theorem~\ref{theosuff}.

\begin{Lemma}[Continuity of solutions]
\label{cont}
Let $u^{\epsilon}$ be a control perturbation around the optimal 
control $u^{*}$, that is, for all $t\in [a,b]$, 
$u^{\epsilon}(t)= u^{*}(t)+\epsilon h(t)$, where 
$h(\cdot) \in L^2\left([a, b], \mathbb{R}\right)$ 
is a variation and $\epsilon \in \mathbb{R}$. Denote by $x^{\epsilon}$ 
its corresponding state trajectory, solution of 
\begin{equation}\label{eq:control}
^A_{C}D^{\alpha, \beta, \gamma}_{a, b}x^{\epsilon}(t)
= f\left(t, x^{\epsilon}(t), u^{\epsilon}(t)\right), 
\quad x^{\epsilon}(a)= x_a.
\end{equation}
If $\displaystyle{ K < \frac{1}{(b-a)^{\alpha-1}M}}$, 
where $K$ is the Lipschitz constant of $f$,
then $x^{\epsilon}$ converges to the optimal 
state trajectory $x^{*}$ when $\epsilon$ tends to zero,
that is, $x^{*}$ is continuous.
\end{Lemma}

\begin{proof}
From Equation \eqref{eq:control}, we have
\[
\left| \, ^A_{C}D^{\alpha, \beta, \gamma}_{a, b} x^{\epsilon}(t) 
- ^A_{C}D^{\alpha, \beta, \gamma}_{a, b} x^{*}(t) \right|
= \left|f\left(t, x^{\epsilon}(t), u^{\epsilon}(t)\right)
- f\left(t, x^{*}(t), u^{*}(t)\right)\right|.
\]
Abbreviating $f^{\epsilon}-f^{*}
= f\left(t, x^{\epsilon}(t), u^{\epsilon}(t)\right)
- f\left(t, x^{*}(t), u^{*}(t)\right)$, it follows, 
by definition of combined operators \eqref{def:combined}, that
\begin{equation}
\label{eq:p:c}
\left| \gamma \left( \, ^A_{C}D^{\alpha, \beta}_{a^{+}} x^{\epsilon}(t) 
- ^A_{C}D^{\alpha, \beta}_{a^{+}} x^{*}(t) \right) + (1-\gamma)
\left( ^A_{C}D^{\alpha, \beta}_{b^{-}} x^{\epsilon}(t) 
- ^A_{C}D^{\alpha, \beta}_{b^{-}} x^{*}(t) \right) 
\right|= \left|f^{\epsilon}-f^{*} \right|.
\end{equation}
Next, since $\gamma \in [0, 1]$, we obtain from \eqref{eq:p:c} 
the two separate inequalities
\begin{equation}
\label{inequal1}
\left| \, ^A_{C}D^{\alpha, \beta}_{a^{+}} x^{\epsilon}(t) 
- ^A_{C}D^{\alpha, \beta}_{a^{+}} x^{*}(t) \right| 
\leq \left|f^{\epsilon}-f^{*} \right|
\end{equation}
and 
\begin{equation}
\label{inequal2}
\left| \, ^A_{C}D^{\alpha, \beta}_{b^{-}} x^{\epsilon}(t) 
- ^A_{C}D^{\alpha, \beta}_{b^{-}} x^{*}(t) \right|  
\leq \left|f^{\epsilon}-f^{*} \right|.
\end{equation}
Therefore, considering \eqref{inequal2}, 
we can deduce the integral relation 
\[
\left| x^{\epsilon}(t)- x^{*}(t)\right| 
\leq {^{\bar{A}}I^{\alpha, \beta}_{b^{-}}}\left( |f^{\epsilon}-f^{*}| \right).
\]
By the Lipschitz property of $f$, we determine that for each $t\in [a, b]$, 
there exists $B_1, B_2 \subset \mathbb{R}$, neighbourhood of $x^{*}(t)$, $u^{*}(t)$, 
respectively, and such that 
\begin{align*}
\left| x^{\epsilon}(t)- x^{*}(t)\right| 
&\leq  {^{\bar{A}}I^{\alpha, \beta}_{b^{-}}}\left( K\left| 
x^{\epsilon}(t)- x^{*}(t)\right| + K\left|\epsilon h(t)\right|\right)\\
&= K|\epsilon |^{\bar{A}}I^{\alpha, \beta}_{b^{-}}\left( |h(t)|\right) 
+ K ^{\bar{A}}I^{\alpha, \beta}_{b^{-}}\left( 
\left| x^{\epsilon}(t)- x^{*}(t)\right|\right). 
\end{align*}
Now, applying Gronwall's inequality (Theorem~\ref{Gronwall}), 
with $\gamma = 0$, we have 
\begin{align*}
\left| x^{\epsilon}(t)- x^{*}(t)\right| 
&\leq K|\epsilon |^{\bar{A}}
I^{\alpha, \beta}_{b^{-}}\left( |h(t)|\right) 
+ \sum_{k=1}^{\infty}K^k\left[ ^{\bar{A}}I^{k\alpha,\beta}_{b^{-}}\left(
K|\epsilon|^{\bar{A}}I^{\alpha, \beta}_{b^{-}}\left( |h(t)|\right) \right)  \right]\\
&= |\epsilon|K \left[ ^{\bar{A}}I^{\alpha, \beta}_{b^{-}}\left( |h(t)|\right) 
+ \sum_{k=1}^{\infty}K^k\left( ^{\bar{A}}I^{(k+1)\alpha, \beta}_{b^{-}}\left(
|h(t)|\right)\right)\right].
\end{align*}
Moreover, using a similar method of reasoning, we may consider \eqref{inequal1} and obtain
\[ 
\left| x^{\epsilon}(t)- x^{*}(t)\right| \leq |\epsilon|
K \left[ ^{\bar{A}}I^{\alpha, \beta}_{a^{+}}\left( |h(t)|\right) 
+ \sum_{k=1}^{\infty}K^k\left( ^{\bar{A}}I^{(k+1)\alpha, 
\beta}_{a^{+}}\left( |h(t)|\right)\right)\right].
\]
Hence, summing altogether, we get
\[
\left| x^{\epsilon}(t)- x^{*}(t)\right| 
\leq \frac{1}{2}|\epsilon |K \Upsilon(t),
\]
where 
$$ 
\Upsilon (t)= \left[ \left(^{\bar{A}}I^{\alpha, \beta}_{a^{+}} 
+ ^{\bar{A}}I^{\alpha, \beta}_{b^{-}}\right)\left( |h(t)|\right)\right]  
+ \left[ \sum_{k=1}^{\infty}K^k\left(^{\bar{A}}I^{\alpha, \beta}_{a^{+}} 
+ ^{\bar{A}}I^{\alpha, \beta}_{b^{-}}\right)\left( |h(t)|\right)\right], 
\quad t\in [a, b]. 
$$
Finally, when $\epsilon \rightarrow 0$, we obtain 
$x^{\epsilon}(t) \rightarrow x^{*}(t)$ for all $t \in [a, b]$. 
\end{proof}

\begin{Corollary}[Differentiability of the perturbed trajectory]
\label{different}
There exists a function $\eta$ defined on $[a,b]$ such that 
\[
x^{\epsilon}(t)= x^{*}(t) + \epsilon \eta(t) + o(\epsilon).
\]
\end{Corollary}

\begin{proof}
Since $f \in C^1$, we have that
\begin{multline*}
f(t, x^{\epsilon}, u^{\epsilon})
= f(t, x^{*}, u^{*}) + (x^{\epsilon}-x^{*})
\frac{\partial f(t, x^{*}, u^{*})}{\partial x} 
+ (u^{\epsilon}-u^{*})\frac{\partial f(t, x^{*}, u^{*})}{\partial u}\\
+ o(|x^{\epsilon} - x^{*}|,|u^{\epsilon}-u^{*}|).
\end{multline*}
Observe that $u^{\epsilon}-u^{*}= \epsilon h(t)$, 
$u^{\epsilon} \rightarrow u^{*}$ when $\epsilon \rightarrow 0$, 
and, by Lemma~\ref{cont}, we have $x^{\epsilon} \rightarrow x^{*}$ 
when $\epsilon \rightarrow 0$. Thus, the residue term can be expressed 
in terms of $\epsilon$ only, that is, the residue is $o(\epsilon)$. 
Therefore, 
\[
^A_{C}D^{\alpha, \beta, \gamma}_{a, b}x^{\epsilon}
= {^A_{C}D}^{\alpha, \beta, \gamma}_{a, b}x^{*} + (x^{\epsilon}
-x^{*})\frac{\partial f(t, x^{*}, u^{*})}{\partial x} 
+ \epsilon h(t)\frac{\partial f(t, x^{*}, u^{*})}{\partial u} 
+ o(\epsilon),
\]
which leads to 
\[
\lim_{\epsilon \rightarrow 0} \left[\frac{^A_{C}D^{\alpha, 
\beta, \gamma}_{a, b}( x^{\epsilon}-x^{*})}{\epsilon} 
-\frac{( x^{\epsilon}-x^{*})}{\epsilon}\frac{\partial 
f(t, x^{*}, u^{*})}{\partial x} - h(t)
\frac{\partial f(t, x^{*}, u^{*})}{\partial u}\right]=0, 
\]
that is, 
\[
^A_{C}D^{\alpha, \beta, \gamma}_{a, b} \left(\lim_{\epsilon \rightarrow 0} 
\frac{  x^{\epsilon}-x^{*}}{\epsilon} \right) 
= \lim_{\epsilon \rightarrow 0} \frac{( x^{\epsilon}-x^{*})}{\epsilon}
\frac{\partial f(t, x^{*}, u^{*})}{\partial x} 
+ h(t)\frac{\partial f(t, x^{*}, u^{*})}{\partial u}. 
\]
Now, it remains to prove the existence of the limit
$\displaystyle{\lim_{\epsilon \rightarrow 0} 
\frac{x^{\epsilon}-x^{*}}{\epsilon}} =: \eta$. It is easy to see 
that the limit $\eta$ exists, as a solution of the following 
fractional differential equation:
\begin{equation*}
\begin{cases}
^A_{C}D^{\alpha, \beta, \gamma}_{a, b} \eta(t)
= \frac{\partial f(t, x^{*}, u^{*})}{\partial x}\eta(t) 
+\frac{\partial f(t, x^{*}, u^{*})}{\partial u}h(t),\\[3mm]
\eta(a)= 0.
\end{cases}
\end{equation*}
This ends the proof.
\end{proof}

The following result is a necessary optimality condition for 
the analytic kernel fractional optimal control problem \eqref{bp}.

\begin{Theorem}[Pontryagin Maximum Principle for \eqref{bp}]
\label{theo}
If $(x^{*}(\cdot), u^{*}(\cdot))$ is an optimal pair for~\eqref{bp}, 
then there exists $\lambda \in L^{\alpha, \beta}\left([a, b], \mathbb{R}\right)$, 
called the adjoint function variable, such that the following conditions 
hold in the interval $[a, b]$:
\begin{itemize}
\item the optimality condition
\begin{equation}\label{opt}
u^{*} \, \, \text{maximizes, over} 
\, \, \Omega(t), \, \, \text{the function}, \, \, 
u \mapsto H\left(t, x^{*}(t), u, \lambda(t) \right);
\end{equation}
\item the adjoint equation
\begin{equation}
\label{adj}
^{A}_{RL}D^{\alpha, \beta, 1-\gamma}_{a, b} \lambda (t)
= \frac{\partial H}{\partial x}(t,x^{*}(t), u^{*}(t)) 
;
\end{equation}
\item the transversality condition
\begin{equation}
\label{trans}
(1-\gamma)^{\bar{A}}I^{1-\alpha, \beta}_{a^{+}}\lambda(b) 
+ \gamma ^{\bar{A}}I^{1-\alpha, \beta}_{b^{-}}\lambda(b)=0,
\end{equation}
where $\bar{A}$ is such that 
$A_{\Gamma}(x^{*}) \cdot \bar{A}_{\Gamma}(x^{*})=1$ and $H(t, x, u, \lambda)= L + \lambda f.$
\end{itemize}
\end{Theorem}

\begin{proof}
Let $(x^{*}(\cdot), u^{*}(\cdot))$ be solution of problem \eqref{bp}, 
and $h\in  L^2\left([a, b], \mathbb{R}\right)$ be a variation, that is, 
\[
h \in V := \left\{v \in L^2: u^{*} + \epsilon v \in \Omega (t), \, \, 
\text{ for any } \epsilon > 0 \, \, \text{ sufficiently small } \right\}.
\]
Set $u^{\epsilon}(t)= u^{*}(t)+ \epsilon h(t)$, so that  
$u^{\epsilon}\in \Omega(t)$ and let $x^{\epsilon}$ be the state 
corresponding to the control $u^{\epsilon}$, that is, the solution of 
\begin{equation}
\label{equatepsi}
^{A}_{C}D^{\alpha, \beta, \gamma}_{a,b} x^{\epsilon}(t)
= f\left(t, x^{\epsilon}(t), u^{\epsilon}(t)\right), 
\quad x^{\epsilon}(a)= x_a.
\end{equation}
Note that $u^{\epsilon}(t) \rightarrow u^{*}(t)$ for all $t\in [a, b]$ 
whenever $\epsilon \rightarrow 0$. Furthermore, 
\begin{equation}
\label{eqpartialu}
\displaystyle{\left.\frac{\partial 
u^{\epsilon}(t)}{\partial \epsilon}\right|_{\epsilon=0}  = h(t)}.
\end{equation} 
Something similar is also true for $x^{\epsilon}$: 
this is justified by Lemma~\ref{cont}. Indeed, 
because $f\in C^1$ with Lipschitz constant $K$ satisfying 
$K < \displaystyle{\frac{1}{(b-a)^{\alpha-1}M}}$, it follows from 
Lemma~\ref{cont} that, for each fixed $t$,  
$x^{\epsilon}(t)\rightarrow x^{*}(t)$ as $ \epsilon \rightarrow 0$. 
Moreover, by Corollary~\ref{different}, the derivative 
$\displaystyle{\left.\frac{\partial x^{\epsilon}(t)}{\partial \epsilon}\right|_{\epsilon=0}}$ 
exists for each $t$. The objective functional at $(x^{\epsilon}, u^{\epsilon})$ is 
\[
J(x^{\epsilon}, u^{\epsilon})
= \int^b_a L\left(t, x^{\epsilon}(t), u^{\epsilon}(t) \right)dt.
\]
Next, we introduce the adjoint function $\lambda$. 
Let $\lambda (\cdot)$ be in $L^{\alpha, \beta}\left([a, b], \mathbb{R}\right)$, 
to be determined. By the integration by parts Formula \eqref{combineIntpart}, 
\begin{multline*}
\int^b_a \lambda(t)\cdot ^{A}_{C}D^{\alpha, \beta, \gamma}_{a, b} x^{\epsilon}(t)dt 
= \gamma\left[ x^{\epsilon}(t)\cdot ^{\bar{A}}I^{1-\alpha, 
\beta}_{b^{-}}\lambda(t) \right]^b_a + (1-\gamma)\left[ x^{\epsilon}(t)
\cdot ^{\bar{A}}I^{1-\alpha, \beta}_{a^{+}}\lambda(t) \right]^b_a\\ 
+ \int^b_ax^{\epsilon}(t)\cdot ^{A}_{RL}D^{\alpha, \beta, 1-\gamma}_{a, b} 
\lambda (t)dt
\end{multline*}
and one has
\begin{multline*}
\int^b_a \lambda(t)\cdot ^{A}_{C}D^{\alpha, \beta, \gamma}_{a, b} x^{\epsilon}(t)dt 
- \gamma\left[ x^{\epsilon}(t)\cdot ^{\bar{A}}I^{1-\alpha, \beta}_{b^{-}}\lambda(t) \right]^b_a 
- (1-\gamma)\left[ x^{\epsilon}(t)\cdot ^{\bar{A}}I^{1-\alpha, \beta}_{a^{+}}\lambda(t) \right]^b_a\\ 
- \int^b_ax^{\epsilon}(t)\cdot ^{A}_{RL}D^{\alpha, \beta, 1-\gamma}_{a, b} \lambda (t)dt=0.
\end{multline*}
Adding this zero to the expression of $J(x^{\epsilon}, u^{\epsilon})$ gives 
\begin{equation*}
\begin{split}
\phi (\epsilon)&= J(x^{\epsilon}, u^{\epsilon})\\
&= \int^b_a \left[L\left( t, x^{\epsilon}(t), u^{\epsilon}(t)\right) 
+ \lambda(t)\cdot ^{A}_{C}D^{\alpha, \beta, \gamma}_{a, b} x^{\epsilon}(t)
- x^{\epsilon}(t)\cdot ^{A}_{RL}D^{\alpha, \beta, 1-\gamma}_{a, b} 
\lambda (t) \right]dt \\
&\quad - x^{\epsilon}(b)\cdot \left[ (1-\gamma)^{\bar{A}}
I^{1-\alpha, \beta}_{a^{+}}\lambda(b) 
+ \gamma ^{\bar{A}}I^{1-\alpha, \beta}_{b^{-}}\lambda(b)\right]\\
&\quad + x^{\epsilon}(a)\cdot \left[ (1-\gamma)^{\bar{A}}I^{1-\alpha, \beta}_{a^{+}}
\lambda(a) + \gamma ^{\bar{A}}I^{1-\alpha, \beta}_{b^{-}}\lambda(a)\right], 
\end{split}
\end{equation*}
which by \eqref{equatepsi} is equivalent to
\begin{equation*}
\begin{split}
\phi (\epsilon)&= J(x^{\epsilon}, u^{\epsilon})\\
&= \int^b_a \left[L\left( t, x^{\epsilon}(t), u^{\epsilon}(t)\right) 
+ \lambda(t)\cdot f\left(t, x^{\epsilon}(t), u^{\epsilon}(t)\right)
- x^{\epsilon}(t)\cdot ^{A}_{RL}D^{\alpha, \beta, 1-\gamma}_{a, b} \lambda (t) \right]dt\\ 
&\quad - x^{\epsilon}(b)\cdot \left[ (1-\gamma)^{\bar{A}}I^{1-\alpha, \beta}_{a^{+}}\lambda(b) 
+ \gamma ^{\bar{A}}I^{1-\alpha, \beta}_{b^{-}}\lambda(b)\right]\\ 
&\quad + x^{\epsilon}(a) \cdot \left[ (1-\gamma)^{\bar{A}}I^{1-\alpha, \beta}_{a^{+}}\lambda(a) 
+ \gamma ^{\bar{A}}I^{1-\alpha, \beta}_{b^{-}}\lambda(a)\right]. 
\end{split}
\end{equation*}
Since the maximum of $J$ occurs at $(x^{*}, u^{*})= (x^0, u^0)$, we have 
that for every feasible direction's variation (i.e., $h\in V$), the 
derivative of $\phi(\epsilon)$ with respect to $\epsilon$ at $\epsilon=0$ 
must be negative~\cite{MR3979978}, that is,

\begin{multline*}
0\geqslant \phi'(0)= \frac{d }{d \epsilon } 
J(x^{\epsilon}, u^{\epsilon})|_{\epsilon=0}
= \int^b_a \left[ \frac{\partial L}{\partial x}
\left.\frac{\partial x^{\epsilon}(t)}{\partial \epsilon}\right|_{\epsilon=0} 
+ \frac{\partial L}{\partial u}\left.\frac{\partial u^{\epsilon}(t)}{\partial 
\epsilon}\right|_{\epsilon=0} \right. \\  
\left. +  \lambda (t)\left( \frac{\partial f}{\partial x}
\left.\frac{\partial x^{\epsilon}(t)}{\partial \epsilon}\right|_{\epsilon=0} 
+ \frac{\partial f}{\partial u}\left.\frac{\partial u^{\epsilon}(t)}{\partial 
\epsilon}\right|_{\epsilon=0}\right)
- ^{A}_{RL}D^{\alpha, \beta, 1-\gamma}_{a, b} 
\lambda (t)\left.\frac{\partial x^{\epsilon}(t)}{\partial \epsilon}\right|_{\epsilon=0} \right]dt \\
- \left[ (1-\gamma)^{\bar{A}}I^{1-\alpha, \beta}_{a^{+}}\lambda(b) + \gamma ^{\bar{A}}I^{1-\alpha,
\beta}_{b^{-}}\lambda(b)\right] \left.
\frac{\partial x^{\epsilon}(b)}{\partial \epsilon}\right|_{\epsilon=0},
\end{multline*}
\textls[-21]{where the partial derivatives of $L$ and $f$ with respect to $x$ and $u$ are evaluated at 
$\left( t, x^{*}(t), u^{*}(t) \right)$. Rearranging the terms and using \eqref{eqpartialu}, 
we obtain that}
\begin{multline*}
\int^b_a \left[ \left(\frac{\partial L}{\partial x} + \lambda(t)\frac{\partial f}{\partial x } 
-  ^{A}_{RL}D^{\alpha, \beta, 1-\gamma}_{a, b} \lambda (t)\right) \left.\frac{\partial 
x^{\epsilon}(t)}{\partial \epsilon}\right|_{\epsilon=0} + \left( \frac{\partial L}{\partial u} 
+ \lambda(t)\frac{\partial f}{\partial u}\right)h(t)\right]dt \\
- \left[ (1-\gamma)^{\bar{A}}I^{1-\alpha, \beta}_{a^{+}}\lambda(b) 
+ \gamma ^{\bar{A}}I^{1-\alpha, \beta}_{b^{-}}\lambda(b)\right] 
\left.\frac{\partial x^{\epsilon}(b)}{\partial \epsilon}\right|_{\epsilon=0}\leqslant 0.
\end{multline*}
Setting $H= L + \lambda f$, it follows that
\begin{multline*}
\int^b_a \left[ \left(\frac{\partial H}{\partial x} 
- ^{A}_{RL}D^{\alpha, \beta, 1-\gamma}_{a, b} \lambda (t)\right) 
\left.\frac{\partial x^{\epsilon}(t)}{\partial \epsilon}\right|_{\epsilon=0} 
+ \frac{\partial H}{\partial u} h(t)\right]dt\\ 
-  \left[ (1-\gamma)^{\bar{A}}I^{1-\alpha, \beta}_{a^{+}}\lambda(b)
+ \gamma ^{\bar{A}}I^{1-\alpha, \beta}_{b^{-}}\lambda(b)\right]
\left.\frac{\partial x^{\epsilon}(b)}{\partial \epsilon}\right|_{\epsilon=0}\leqslant 0,
\end{multline*}
where the partial derivatives of $H$ are evaluated at
$\left( t, x^{*}(t), u^{*}(t), \lambda(t) \right)$. 
Now, choosing
\[
^{A}_{RL}D^{\alpha, \beta, 1-\gamma}_{a, b} \lambda (t) 
= \frac{\partial H}{\partial x}, \quad \text{ with } 
\left[ (1-\gamma)^{\bar{A}}I^{1-\alpha, \beta}_{a^{+}}\lambda(b) 
+ \gamma ^{\bar{A}}I^{1-\alpha, \beta}_{b^{-}}\lambda(b)\right]=0,
\]
that is, given the adjoint Equation \eqref{adj} 
and the transversality condition \eqref{trans}, 
it yields 
\[
\int^b_a \frac{\partial H}{\partial 
u}\left(t, x^{*}(t), u^{*}(t), \lambda(t) \right) h(t) dt \leqslant 0.
\]
Since this inequality holds for any feasible direction's variation 
$h(\cdot) \in V$, we obtain that the partial derivative 
$\frac{\partial H}{\partial u}\left(t, x^{*}(t), u^{*}(t),  
\lambda(t)\right)$ belongs to the normal cone to $\Omega(t)$ at $u^{*}$ 
(see, e.g., p. 45 of~\cite{barbu}), that is, mathematically,
\[
\frac{\partial H}{\partial u}\left(t, x^{*}(t), 
u^{*}(t), \lambda(t)\right) \in N_{\Omega(t)}(u^{*}),
\]
meaning that $u^{*}$ maximizes $u \mapsto H\left(t, x^{*}(t), 
u, \lambda(t)\right)$ over $\Omega(t)$, which is exactly 
the optimality condition \eqref{opt}.
This completes the proof.
\end{proof}

\begin{Remark}
\label{rem:1}
If $\Omega(t)= L^2$, then $ N_{\Omega(t)}(u^{*})={0}$, 
and the optimality condition \eqref{opt} is reduced to 
\[
\frac{\partial H}{\partial u}\left(t, 
x^{*}(t), u^{*}(t), \lambda(t)\right) =0,
\]
which gives the particular result obtained in~\cite{Ndairou:Torres:2021}
(see also \cite{faical}, in a different context).
\end{Remark}

\begin{Definition}
The candidates to solutions of \eqref{bp}, obtained by 
the application of our Theorem~\ref{theo}, 
will be called (Pontryagin) extremals.
\end{Definition}

\begin{Example}
\label{Pexemple}
Let us consider the following optimal control problem:
\begin{equation}
\label{pex}
\begin{gathered}
J[x(\cdot), u(\cdot)]= \int^{2}_{0} -\left(x(t)-t^3 \right)^2 
-\left(u(t) - \frac{t(t-2)}{\ln t} \right)^2 dt \longrightarrow \max,\\
^A_{C}D^{\frac{1}{2}, 3i, \gamma}_{0, 2}x(t)= \frac{1}{2}x(t) 
+ \frac{1}{2}u(t), \quad t\in [0, 2],\\
x(\cdot) \in AC^{\frac{1}{2}, 3i}\left([0, 2], \mathbb{R}\right), 
\quad u(\cdot) \in L^2\left([a, b], \mathbb{R}\right),\\
x(0)= 0,
\end{gathered}
\end{equation}
with 
\[
A(x)=\sum_{n=0}^{\infty}\frac{x^n}{\exp(\sqrt{2})n !}
=\exp(x-\sqrt{2}), \quad \text{ for which }
\quad  M = \underset{|x|<\sqrt{2}}\sup \exp(x-\sqrt{2})=1.
\]
Identifying \eqref{pex} with \eqref{bp}, we have the following correspondence: 
$a=0$; $b=2$; $\displaystyle{\alpha = \frac{1}{2}}$; $\beta = 3i$;  
$f\left(t, x, u \right)=\displaystyle{\frac{1}{2}x + \frac{1}{2}u}$; 
and $\displaystyle{L\left(t, x, u \right) 
= -\left(x-t^3 \right)^2 -\left(u- \frac{t(t-2)}{\ln t} \right)^2}$. 
Note that $f \in C^1$ is Lipschitz-continuous in both variables $x$ and $u$ 
with Lipschitz constant $\displaystyle{K= \frac{1}{2}}$. Also, the inequality 
$\displaystyle{K < \frac{1}{(b-a)^{\alpha-1}M}}$ holds, that is, 
$\displaystyle{\frac{1}{2}< \sqrt{2}}$
is satisfied. Thus, by Lemma~\ref{cont}, for any control perturbation 
the corresponding state trajectory converges to the optimal state solution.
Moreover, by defining the Hamiltonian function as
\begin{equation}
\label{hami}
H(t, x, u,  \lambda) = -\left(x-t^3 \right)^2 
-\left(u- \frac{t(t-2)}{\ln t} \right)^2 + \frac{1}{2}\lambda (x + u),
\end{equation}
and applying Theorem~\ref{theo}, it follows:
\begin{itemize}
\item from the  optimality condition 
$\displaystyle{ \frac{\partial H}{\partial u}= 0}$
(recall Remark~\ref{rem:1}), that
\begin{equation}
\label{optex}
\lambda (t)= 4\left(u(t) - \frac{t(t-2)}{\ln t} \right);
\end{equation}
\item from the adjoint equation 
$\displaystyle{^{A}_{RL}D^{\frac{1}{2}, 3i, 1-\gamma}_{0, 2} \lambda (t)
= \frac{\partial H}{\partial x}}$, that
\begin{equation}
\label{adjex}
^{A}_{RL}D^{\frac{1}{2}, 3i,1-\gamma}_{0, 2} \lambda (t)
=-2\left(x(t)-t^3 \right) + \frac{1}{2}\lambda(t) ;
\end{equation}
\item  from the transversality condition
\begin{equation*}
(1-\gamma)^{\bar{A}}I^{1-\alpha, \beta}_{a^{+}}\lambda(b) 
+ \gamma ^{\bar{A}}I^{1-\alpha, \beta}_{b^{-}}\lambda(b)=0,
\end{equation*}
that 
\begin{equation}
\label{transex}
\left((1-\gamma)^{\bar{A}}I^{\frac{1}{2}, 3i}_{0^{+}} 
+ \gamma ^{\bar{A}}I^{\frac{1}{2}, 3i}_{2^{-}}\right)\lambda(2)=0.
\end{equation}
\end{itemize}
In conclusion, we easily see that \eqref{optex}--\eqref{transex} 
are satisfied by the triple
\begin{equation}
\label{eq:PE}
x(t)= t^3, \quad u(t)= \frac{t(t-2)}{\ln t}, 
\quad \text{ and } \quad \lambda(t) \equiv 0,
\end{equation}
which is the Pontryagin extremal: a candidate 
to the solution of the given problem \eqref{pex}.
\end{Example}


\subsection{Sufficient Condition for Global Optimality}
\label{subsec:SCGO}

We now prove a Mangasarian type theorem for the general 
analytic kernel fractional optimal control problem \eqref{bp}.

\begin{Theorem}[Sufficient global optimality condition]
\label{theosuff}
Consider the general analytic kernel  
fractional optimal control problem \eqref{bp}
with $\Omega(t)= L^2$.
If $(x, u) \rightarrow L(t, x, u)$ and $(x, u) \rightarrow f(t, x, u)$
are concave and $(\tilde{x}, \tilde{u}, \lambda)$ is a Pontryagin extremal 
with $\lambda (t)\geq 0$, $t \in [a,b]$, then  
$$
J[\tilde{x}, \tilde{u}] \geq J[x,u] 
$$ 
for any admissible pair $(x,u)$, that is, 
the pair $\left(\tilde{x}, \tilde{u}\right)$ 
is the solution to problem \eqref{bp}.
\end{Theorem}

\begin{proof}
Since $L$ is concave as a function of $x$ and $u$, we have 
by the gradient inequality (Lemma~\ref{lemma:concave}) that
\begin{multline*}
L\left(t, \tilde{x}(t), \tilde{u}(t)\right)
- L\left(t, x(t), u(t)\right)
\geq \frac{\partial L}{\partial x}\left(t, \tilde{x}(t), 
\tilde{u}(t)\right)\cdot \left(\tilde{x}(t)- x(t)\right)\\ 
+ \frac{\partial L}{\partial u}\left(t, \tilde{x}(t), 
\tilde{u}(t)\right)\cdot \left(\tilde{u}(t)- u(t)\right)
\end{multline*}
for any control $u$ and its associated trajectory $x$. 
This gives 
\begin{equation}
\label{difrfunctio}
\begin{split}
& J[\tilde{x}(\cdot), \tilde{u}(\cdot)]- J[x(\cdot), u(\cdot)]
= \int^{b}_{a}\left[L\left(t, \tilde{x}(t), \tilde{u}(t)\right)
- L\left(t, x(t), u(t)\right)\right]dt\\ 
& \geq \int^b_a \left[ \frac{\partial L}{\partial x}
\left(t, \tilde{x}(t), \tilde{u}(t)\right)\cdot \left(\tilde{x}(t)- x(t)\right) 
+ \frac{\partial L}{\partial u}\left(t, \tilde{x}(t), 
\tilde{u}(t)\right)\cdot \left(\tilde{u}(t)- u(t)\right)\right]dt\\ 
&= \int^b_a \left[\frac{\partial \tilde{L}}{\partial x}
\cdot \left(\tilde{x}(t)-x(t) \right) 
+ \frac{\partial \tilde{L}}{\partial u}
\cdot \left(\tilde{u}(t)-u(t) \right) \right]dt,
\end{split}
\end{equation}
where $\tilde{L}= L\left(t, \tilde{x}(t), \tilde{u}(t)\right)$.
From the adjoint Equation \eqref{adj}, we can write
\[
\frac{\partial L}{\partial x}(t, \tilde{x}(t), \tilde{u}(t)) 
= {^{A}_{RL}D^{\alpha, \beta, 1-\gamma}_{a, b}} \lambda(t) 
- \lambda(t)\frac{\partial f}{\partial x}(t, \tilde{x}(t), \tilde{u}(t)),
\]
while from the optimality condition \eqref{opt} (recall Remark~\ref{rem:1}) we have  
\[
\frac{\partial L}{\partial u}(t, \tilde{x}(t), \tilde{u}(t))
= - \lambda (t)\frac{\partial f}{\partial u}(t, \tilde{x}(t), \tilde{u}(t)).
\]
It follows from \eqref{difrfunctio} that
\begin{multline}
\label{difrfunctio1}
J[\tilde{x}(\cdot),  \tilde{u}(t)]- J[x(\cdot ), u(\cdot)]\\
\geq \int^b_a \left[ \left( ^{A}_{RL}D^{\alpha, \beta, 1-\gamma}_{a, b}\lambda (t) 
- \lambda (t)\frac{\partial \tilde{f}}{\partial x}\right)
\cdot \left( \tilde{x}(t)-x(t) \right) - \lambda (t)
\frac{\partial \tilde{f}}{\partial u} \cdot \left( \tilde{u}(t) -u(t)\right) \right]dt,
\end{multline} 
where $\tilde{f}= f\left(t, \tilde{x}(t), \tilde{u}(t)\right)$. 
Next, by using the integration by parts Formula \eqref{combineIntpart}, we get
\begin{multline*}
\int^b_a \lambda(t)\cdot ^{A}_{C}D^{\alpha, \beta, \gamma}_{a, b}\left( 
\tilde{x}(t) -x(t)\right)dt = \gamma\left[ \left( \tilde{x}(t) -x(t)\right)
\cdot {^{\bar{A}}I^{1-\alpha, \beta}_{b^{-}}} \lambda(t) \right]^b_a \\ 
+ (1-\gamma)\left[ \left( \tilde{x}(t) -x(t)\right)
\cdot {^{\bar{A}} I^{1-\alpha, \beta}_{a^{+}}} \lambda(t) \right]^b_a  
+ \int^b_a \left( \tilde{x}(t) -x(t)\right)
\cdot ^{A}_{RL}D^{\alpha, \beta, 1-\gamma}_{a,b}\lambda (t)dt,
\end{multline*}
meaning that 
\begin{multline}
\label{integra}
\int^b_a \left( \tilde{x}(t) -x(t)\right)\cdot ^{A}_{RL}
D^{\alpha, \beta, 1-\gamma}_{a,b}\lambda (t)dt
= \int^b_a \lambda(t)\cdot ^{A}_{C}D^{\alpha, \beta, \gamma}_{a, b}\left( 
\tilde{x}(t) -x(t)\right)dt\\ - \gamma\left[ \left( \tilde{x}(t) 
-x(t)\right)\cdot ^{\bar{A}}I^{1-\alpha, \beta}_{b^{-}}\lambda(t) \right]^b_a 
- (1-\gamma)\left[ \left( \tilde{x}(t) -x(t)\right)\cdot ^{\bar{A}}
I^{1-\alpha, \beta}_{a^{+}}\lambda(t) \right]^b_a.  
\end{multline}
Substituting \eqref{integra} into \eqref{difrfunctio1}, we get
\begin{equation*}
J\left[ \tilde{x}(\cdot), \tilde{u}(\cdot) \right] 
- J\left[ x(\cdot), u(\cdot)\right]
\geq \int^b_a \lambda(t)\left[ \tilde{f} 
- f - \frac{\partial \tilde{f}}{\partial x}\left( \tilde{x}(t)- x(t)\right) 
- \frac{\partial \tilde{f}}{\partial u} \left( \tilde{u}(t)- u(t)\right)\right]dt.
\end{equation*}
Finally, taking into account that $\lambda(t)\geq 0$ and $f$ is concave 
in both $x$ and $u$, we conclude that 
$J\left[ \tilde{x}(\cdot), \tilde{u}(\cdot) \right] 
- J\left[ x(\cdot), u(\cdot)\right]\geq 0$.
\end{proof}

\begin{Example}
\label{suff:example}
It is easily proved from Theorem~\ref{theosuff} 
that the Pontryagin extremal \eqref{eq:PE},
$$
x(t)= t^3, \quad u(t)= \frac{t(t-2)}{\ln t}, 
\quad \text{ and } \quad \lambda(t)\equiv 0,
$$  
candidate to the solution of the optimal control problem \eqref{pex},
found in Example~\ref{Pexemple} from the application of Theorem~\ref{theo},
is indeed a solution to the problem (it is a global maximizer):  
in this case, the Hamiltonian defined in \eqref{hami} is a concave function 
with respect to both variables $x$ and $u$ and, furthermore, 
$\lambda (t)\geq 0$ for all $t \in [a,b]$.
\end{Example}


\section{Conclusions}
\label{sec:conc}

In this paper, we investigated, for the first time in the literature,
optimal control problems with combined general analytic kernels.
Main results provide strong necessary optimality conditions
of Pontryagin type, valid in the class of 
absolutely continuous state trajectories
and $L^2$ controls that may take values
in any time-dependent close convex set.
Other results include a new Gronwall inequality 
and a sufficient optimality condition for global maximizers.
While our results provide non-trivial and useful analytical
results, as here shown with a simple illustrative example, 
to address real-world applications it will be necessary
to develop numerical methods that implement the obtained results.
This opens several possible future directions of research
and will be addressed elsewhere.


\vspace{6pt} 


\authorcontributions{The authors equally contributed to this paper, read and approved the final
manuscript. Formal analysis, F.N. and D.F.M.T.; Investigation, F.N. and D.F.M.T.; Writing---original
draft, F.N. and D.F.M.T.; Writing---review \& editing, F.N. and D.F.M.T. Both authors have read and
agreed to the published version of the manuscript.}
	
\funding{This research was funded by the Portuguese Foundation for Science and Technology (FCT),
grant number UIDB/04106/2020 (CIDMA). Nda\"{\i}rou was also supported by FCT through the PhD
fellowship PD/BD/150273/2019.}

\dataavailability{Not applicable.} 

\acknowledgments{The authors are grateful to two anonymous reviewers
for several constructive remarks and questions.}

\conflictsofinterest{The authors declare no conflict of interest.
The funders had no role in the design of the study; in the collection, 
analyses, or interpretation of data; in the writing of the manuscript, 
or in the decision to publish the~results.}

\end{paracol}


\reftitle{References}


\end{document}